\documentstyle{amsppt}
\magnification=1200
\hsize=150truemm
\vsize=224.4truemm
\hoffset=4.8truemm
\voffset=12truemm

\NoRunningHeads

\newsymbol\lessim 132E
\define\C{{\bold C}}
\let\th\proclaim
\let\fth\endproclaim
 
\let\bo\partial

\newcount\refno
\global\refno=0

\def\nextref#1{
      \global\advance\refno by 1
      \xdef#1{\the\refno}}

\def\bref {\ref\global\advance\refno by 1\key{\the\refno}}


\nextref\BE
\nextref\CH
\nextref\E
\nextref\ES
\nextref\G
\nextref\H
\nextref\K
\nextref\M
\nextref\N
\nextref\OK
\nextref\R

\topmatter
\title 
Sur les courbes de Brody dans $P^n(\C)$
 \endtitle
\author Bernardo Freitas Paulo Da Costa et Julien Duval 
\footnote{Laboratoire de Math\'ematiques, Universit\'e Paris-Sud, 91405 Orsay cedex, France \newline bernardo.da-costa\@math.u-psud.fr\newline julien.duval\@math.u-psud.fr}
\footnote""{Mots cl\'es : courbes de Brody, th\'eorie de Nevanlinna  \newline MSC 2010: 30D35, 32H25 \newline}
\endauthor
 
\abstract We investigate Brody curves in the projective space from the point of view of Nevanlinna theory. In particular we show that Brody curves in $P^1(\C)$ have no defect whenever of order $>1$.
\endabstract
  \endtopmatter 
\document
 \subhead 0. Introduction \endsubhead

\null
Ce travail porte sur la th\'eorie de Nevanlinna des courbes de Brody dans $P^n(\C)$. Une courbe de Brody dans une vari\'et\'e compacte complexe $X$ est une courbe enti\`ere non constante de d\'eriv\'ee born\'ee. L'int\'er\^et de ces courbes provient du lemme de Brody (voir [\K]). Il affirme que toute suite divergente de disques holomorphes dans $X$ produit par renormalisation une courbe de Brody. Ainsi toute vari\'et\'e compacte complexe non hyperbolique au sens de Kobayashi contient de telles courbes.

Commen\c cons par $P^1(\C)$. Soit $f :\C \to P^1(\C)$ une courbe enti\` ere. La th\'eorie de Nevanlinna [\N] compare la croissance de $f$ avec la r\'epartition de ses pr\'eimages. La croissance de $f$ se mesure par l'indicatrice $$T(r)=\int_0^r \left(\int_{D_t} f^*\omega \right)\, \frac{dt}t$$ o\`u $D_t$ est le disque de centre $0$ et de rayon $t$ dans $\C$
et $\omega$ la forme de Fubini-Study normalis\'ee de $P^1(\C)$. L'ordre de $f$ est $\limsup \frac {\log T}{\log r}$. Si $a$ est dans $P^1(\C)$, la r\'epartition de ses pr\'eimages se traduit par la fonction de comptage (avec multiplicit\'es) $$N_a(r)=\int_0^r \text{card}(f^{-1}(a)\cap D_t)\, \frac {dt}t.$$ Quand $f$ est une courbe de Brody, $f^*\omega$ reste born\'ee et $T=O(r^2)$. En particulier $f$ est d'ordre au plus 2. Notre point de d\'epart est le r\'esultat de Clunie et Hayman [\CH] montrant qu'une courbe de Brody de $P^1(\C)$ \'evitant un point a une croissance plus faible : $T=O(r)$. Elle est d'ordre au plus 1. On peut le voir comme un analogue du th\'eor\`eme de Picard qui dit qu'une courbe enti\`ere \'evitant 3 points est constante. Il est alors naturel d'esp\'erer une version pour les courbes de Brody du second th\'eor\`eme fondamental de la th\'eorie de Nevanlinna. Celui-ci quantifie Picard en affirmant que les impacts d'une courbe enti\`ere sur $q$ points distincts croissent au moins comme $(q-2)T$. Voici notre r\'esultat pour les courbes de Brody.

\th{ Th\'eor\`eme 1} Soit $f$ une courbe de Brody dans $P^1(\C)$.
Alors $T=N_a+O(r)$ pour tout $a$ dans $P^1(\C)$. \fth

Donc si $f$ cro\^\i t suffisamment elle est r\'eguli\`ere au sens o\`u la r\'epartition de ses pr\'eimages suit son indicatrice de croissance. Elle se rapproche des prototypes de courbes de Brody d'ordre 2 que sont les fonctions elliptiques.
Pr\'ecisons ceci. Soit $\delta_a= \liminf (1-\frac {N_a} {T})$ le d\'efaut de $f$ en $a$. Le second th\'eor\`eme principal de Nevanlinna donne une majoration par 2 de la somme des d\'efauts. Dans le cas des courbes de Brody d'ordre $>1$ on a mieux.

\th{Corollaire} Soit $f$ une courbe de Brody dans $P^1(\C)$ d'ordre $>1$. Alors elle est sans d\'efaut, $\delta_a=0$ pour tout $a$ dans $P^1(\C)$.\fth

Passons \`a la dimension sup\'erieure. Soit $f:\C \to P^n(\C)$ une courbe enti\`ere. On d\'efinit comme plus haut son indicatrice de croissance $T$ et la fonction de comptage $N_H$ des impacts de $f$ sur un hyperplan $H$ ne contenant pas $f(\C)$. Une courbe de Brody est encore d'ordre au plus 2 mais pour les courbes de Brody \'evitant $n$ hyperplans en position g\'en\'erale Eremenko [\E] a montr\'e l'analogue du r\'esultat de Clunie et Hayman : elles sont d'ordre au plus 1. De plus il y a beaucoup de courbes de Brody d'ordre $>1$ (m\^eme 2) \'evitant $n-1$ hyperplans en position g\'en\'erale. On peut s'int\'eresser \`a leur distribution des valeurs.  Voici notre r\'esultat en dimension sup\'erieure.

\th{ Th\'eor\`eme 2} Soit $f :\C \to P^n(\C)$ une courbe de Brody \'evitant $n-1$ hyperplans en position g\'en\'erale. Soit $H$ un hyperplan en position g\'en\'erale par rapport aux pr\'ec\'edents et ne contenant pas $f(\C)$.
Alors $T=N_H+O(r)$. \fth

Par exemple, soit $f$ une courbe de Brody de $P^2(\C)$ d'ordre $>1$ \'evitant une droite. Alors elle est sans d\'efaut sur toutes les autres droites ne la contenant pas.

\null Pour une courbe de Brody quelconque dans $P^n(\C)$ on s'attend \`a un second th\'eor\`eme de la forme $(q-n+1)\, T\leq \sum N_{H_i}+O(r)$ si les $H_i$ sont $q$ hyperplans en position g\'en\'erale. Nous obtenons dans cette direction un r\'esultat plus faible avec un reste en $o(r^2)$.

\th{ Th\'eor\`eme 3} Soit $f :\C \to P^n(\C) $ une courbe de Brody et $H_i$ $q$ hyperplans en position g\'en\'erale ne contenant pas $f(\C)$. Alors  $(q-n+1)\, T\leq \sum N_{H_i}+o(r^2)$ \fth 

Il implique la formule suivante des d\'efauts pour les courbes de Brody de croissance maximale.

\th{Corollaire} Soit $f$ une courbe de Brody dans $P^n(\C)$ de croissance maximale ($\limsup \frac {T}{r^2}>0$) et $H_i$ des hyperplans en position g\'en\'erale ne contenant pas $f(\C)$. Alors $\sum \delta_{H_i}\leq n-1$.\fth

Notons que la th\'eorie de Cartan (voir [\K]) donne une somme des d\'efauts major\'ee par $n+1$ pour des courbes arbitraires (non d\'eg\'en\'er\'ees) de $P^n(\C) $. Comme dans $P^1(\C)$ on obtient un gain de 2 pour les courbes de Brody.

La d\'emonstration de ces r\'esultats repose sur une r\'eduction \`a des lemmes de potentiel et s'inspire des travaux d'Eremenko et de ses collaborateurs [\E],[\BE],[\ES].

\null\noindent
\subhead 1. Les th\'eor\`emes 1 et 2 \endsubhead

\null\noindent a) {\bf Lemmes de potentiel.} 

\null\noindent

Dans ce paragraphe la m\^eme lettre $C$ peut d\'esigner des constantes universelles diff\'erentes. Si $d$ est un disque du plan de rayon $\rho$, $\lambda d$ est le disque concentrique de rayon $\lambda \rho$. Les int\'egrales sur des arcs de cercle s'entendent pour la mesure de Lebesgue normalis\'ee du cercle. Pour faire court, un {\it arc  de $d$ de rayon $\geq \rho$} est la trace sur $d$ d'un cercle du plan de rayon $\geq \rho$. Les fonctions surharmoniques sont suppos\'es continues (avec \'eventuellement des p\^oles) jusqu'au bord des disques. Notre r\'ef\'erence de base pour le potentiel est [\R].

Commen\c cons par une estim\'ee uniforme due \`a Eremenko [\E] et proche du lemme de Hopf.
\th{Lemme 1} Soit $u$ une fonction surharmonique $\geq 0$ sur $2d$, o\`u $d$ est un disque de rayon $\rho$. On suppose $u$ harmonique sur $\frac{3d}2$ et nulle en un point $z$ au bord de $2d$ o\`u elle est diff\'erentiable. Alors $\text{Max}_d \, u\leq C\rho\,  \vert \nabla u(z)\vert $.\fth

\noindent{\it D\'emonstration.} On se ram\`ene par homoth\'etie et translation au disque  $2D$, o\`u $D$ est le disque unit\'e. Soit $m=\text{Min}_D\, u$ et $M=\text{Max}_D\,  u$. Par le principe du minimum on a $u\geq m e$ o\`u $e$ est la fonction extr\'emale relative de $D$ par rapport \`a $2D$ (valant $1$ sur $D$, 0 au bord de $2D$ et harmonique entre).  Comme $u$ et $e$ s'annulent en $z$ on a $\vert \nabla u(z)\vert \geq m \, \vert \nabla e (z)\vert$. Ceci conclut puisque $M$ et $m$ sont comparables par les in\'egalit\'es de Harnack.

\null 
Sans hypoth\`ese d'harmonicit\'e on peut n\'eanmoins estimer l'int\'egrale de $u$ sur des arcs de cercle assez grands.

\th{Lemme 2} Soit $u$ une fonction surharmonique $\geq0$ sur $2d$, o\`u $d$ est un disque de rayon $\rho$. On suppose $u$ nulle en un point $z$ au bord de $2d$ o\`u elle est diff\'erentiable. On se donne un arc $\gamma$ de $d$ de rayon $\geq  \rho$. Alors $\oint_\gamma u \leq C\rho \, \vert \nabla u(z)\vert$.\fth

\noindent{\it D\'emonstration.} On se ram\`ene \`a $2D$. Soit $g$ la solution du probl\`eme de Dirichlet nulle au bord de $2D$ avec $\Delta g = \Delta u$ sur $\frac{3D}2$. Autrement dit $g=-\int_{\frac{3D}2}g_t \, \Delta u(t)$ o\`u $g_t$ est la fonction de Green surharmonique de $2D$ avec p\^ole en $t$. Par construction $u-g$ est $\geq0$ au bord de $2D$, surharmonique, donc $\geq0$ sur $2D$ par le principe du minimum. Comme elle s'annule en $z$ on a $\vert \nabla g(z)\vert \leq \vert \nabla u(z)\vert$. Par ailleurs $u-g$ est harmonique sur $\frac{3D}2$. Par le lemme 1 on a $u\leq g +C\, \vert \nabla u (z)\vert $ sur $D$. Donc $\oint_\gamma u \leq \int_{\frac{3D}2}(\oint_\gamma g_t)(-\Delta u)+C\, \vert \nabla u (z)\vert $ en int\'egrant. On va estimer l'int\'egrale $I$ de droite. Remarquons d'abord que $\oint_\gamma g_t$ d\'epend contin\^ument de $t$ et du cercle. Elle tend de plus vers z\'ero quand le rayon du cercle tend vers l'infini (on calcule des moyennes). On en d\'eduit par compacit\'e que $\oint_\gamma g_t \leq  C$. Par ailleurs $\vert \nabla g_t(z)\vert \geq \frac{1}{ C}$ quand $t$ parcourt $\frac{3D}2$. Donc
$I\leq C\int_{\frac{3D}2} \vert \nabla g_t(z)\vert\, (-\Delta u)=C\, \vert \nabla g(z)\vert \leq C\, \vert \nabla u(z)\vert $ ce qui conclut.

\null 
Dans ce qui pr\'ec\`ede on peut remplacer le gradient par une int\'egrale sur un arc de cercle centr\'e au bord. On obtient ainsi un r\'esultat de {\it transfert} d'estim\'ees sur des arcs de cercles que le lecteur pourra v\'erifier de la m\^eme mani\` ere. Il sera utile pour la dimension sup\'erieure.

\th{Lemme 3} Soit $u$ une fonction surharmonique $\geq0$ sur $2d$, o\`u $d$ est un disque de rayon $\rho$, et $d'$ un disque centr\'e en un point du bord de $2d$ de rayon $\rho'\leq2\rho$. On pose $\gamma'=2d\, \cap\,  \bo d'$. Alors pour tout arc $\gamma$ de $d$ de rayon $\geq   \rho$ on a $\oint_\gamma u \leq C\frac {\rho}{\rho'} \oint_{\gamma'}u$. \fth

\null 
Toujours pour la dimension sup\'erieure il est important d'avoir un r\'esultat d'{\it oubli} d'une fonction harmonique dans nos estim\'ees sur des arcs de cercles. 
\th{Lemme 4} Soient $u$ et $h$ deux fonctions $\geq0$ respectivement surharmonique et harmonique sur $2d$, o\`u $d$ est un disque de rayon $\rho$. On suppose que $u\leq h$ au centre de $d$. On note $m = \sup_{\Gamma} \, \oint _\gamma \text {Min}(u,h)$ et $M=\sup_{\tilde \Gamma} \, \oint _\gamma u$, o\`u $\Gamma$ (resp. $\tilde \Gamma$) est l'ensemble des arcs de $d$ de rayon $\geq \rho$ (resp. $\geq \frac{\rho}2$). Alors $M\leq C\, m$.\fth

\noindent{\it D\'emonstration.} On se ram\`ene \`a $2D$. Le point cl\'e est le fait suivant.

\null\noindent{\bf Fait.} Soit $v$ surharmonique $\geq0$ sur $2D$. Alors  $\text{Min}_D\, v\leq \sup_\Gamma\, \oint _\gamma v$ et $\sup_{\tilde \Gamma}\, \oint _\gamma v\leq C\, \text{Min}_D\, v$. 
 
\null Admettons ce fait pour l'instant. Voici comment on conclut. En prenant $v= \text {Min}(u,h)$ dans la premi\`ere in\'egalit\'e on obtient un point $z$ de $D$ avec $ \text {Min}(u,h)(z)\leq m$, soit $u(z)\leq m$ ou $h(z)\leq m$. Dans le premier cas on obtient directement l'estim\'ee voulue en posant $v=u$ dans la deuxi\`eme in\'egalit\'e. Dans le second cas on a $h(0)\leq C\, m$ par les in\'egalit\'es de Harnack. Comme par hypoth\`ese $h(0)\geq u(0)$ on conclut de la m\^eme mani\`ere.

\null

Montrons maintenant ce fait. La premi\`ere in\'egalit\'e s'obtient en prenant $\gamma=\bo D$. Pour la deuxi\`eme on suit la d\'emarche du lemme 2. On d\'efinit $g$ de la m\^eme mani\`ere. On suppose le minimum de $v$ sur $D$ atteint en $z$. Par les in\'egalit\'es de Harnack on a $v\leq g+Cv(z)$ sur $D$. Donc par int\'egration $\oint_\gamma v \leq \int_{\frac{3D}2}\, (\oint_\gamma g_t)(-\Delta v)+C\, v(z)$. D'o\`u
$\oint_\gamma v \leq C[-\Delta v(\frac{3D}2)+v(z)]$ puisque les int\'egrales $\oint_\gamma g_t $ sont born\'ees pour $\gamma$ dans $\tilde \Gamma$. Reste \`a estimer la masse de $-\Delta v$ par $v(z)$. Soit $\phi$ un automorphisme de $2D$ envoyant $0$ sur $z$. Notons que $\phi^{-1}(\frac{3D}2)$ reste dans un disque fixe $D'$ plus petit que $2D$. Il suffit d'estimer la masse du laplacien de $w=v\circ \phi$ sur $D'$ par sa valeur au centre. Soit $\chi$ une fonction plateau valant 1 sur $D'$ et \`a support compact dans $2D$. On a $-\Delta w( D')\leq -\int_D \chi\,  \Delta w  =-\int_D w\, \Delta \chi \, dxdy \leq C\int_Dw\, dxdy \leq Cw(0)$ par l'in\'egalit\'e de la moyenne.

\null\noindent
b) {\bf Le th\'eor\`eme 1.} 

\null 
On suit l'approche d'Eremenko [\E] pour le th\'eor\`eme de Clunie et Hayman. Soit $f$ une courbe de Brody
dans $P^1(\C)$. Elle s'\'ecrit $(f_1:f_2)$ en coordonn\'ees homog\`enes avec $f_i$ enti\`eres sans z\'ero commun. On pose $u_i=\log \vert f_i\vert $. L'indicatrice de la courbe est $T=\oint_r \text{Max}(u_1,u_2)$ (\`a un $O(1)$ pr\`es). Ici l'indice $r$ d\'esigne la moyenne sur le cercle $\gamma_r$ centr\'e \`a l'origine de rayon $r$. Le caract\`ere Brody de la courbe se traduit par le fait que $\vert \nabla u_1-\nabla u_2\vert$ est born\'e sur le lieu $(u_1=u_2)$ [\E]. Celui-ci n'est pas vide si $f$ n'est pas constante. On se donne maintenant un point de $P^1(\C)$, par exemple $(1:0)$. Sa fonction de comptage est $N=\oint_r u_2$. Comparons ces deux int\'egrales. On va le faire localement gr\^ace au lemme 2. Fixons un point de $\gamma_r$ o\`u $u_1>u_2$. Soit $2d$ le disque maximal centr\'e en ce point contenu dans $(u_1>u_2)$. Son bord touche $(u_1=u_2)$. Par le lemme 2 appliqu\'e \`a $u_1-u_2$, on en d\'eduit que $\oint_{\gamma}(u_1-u_2) = O(\rho)$ o\`u $\gamma$ d\'esigne la partie de $\gamma_r$ dans $d$ et $\rho$ le rayon de $d$ (contr\^ol\'e par $r$). Par un lemme de recouvrement classique (voir [\M]) on peut trouver un nombre d\'enombrable de tels disques $d_k$ recouvrant $\gamma_r \cap (u_1>u_2)$ tels que les $\frac{d_k}5$ soient disjoints. Il s'ensuit que $\sum \rho_k = O(r)$. En sommant les estim\'ees locales on a donc
$\oint_r(u_1-u_2)^+ =O(r)$ soit $\oint_r\text{Max}(u_1,u_2)=\oint _r u_2+O(r)$.

\null\noindent
c) {\bf Le th\'eor\`eme 2.}

\null  Soit $f$ une courbe de Brody dans $P^n(\C)$. Comme plus haut on peut l'\'ecrire $(f_1:...:f_{n+1})$ en coordonn\'ees homog\`enes avec $f_i$ enti\`eres sans z\'ero commun. On pose $u_i=\log \vert f_i\vert $. L'indicatrice de la courbe est $T=\oint_r \text{Max}(u_1,...,u_{n+1})$. Son caract\`ere Brody se traduit par le fait que $\vert \nabla u_i-\nabla u_j\vert$ est born\'e sur le lieu $(u_i=u_j=\text{Max}(u_1,...,u_{n+1}))$ [\E]. On suppose la courbe non d\'eg\'en\'er\'ee, i.e. non contenue dans un hyperplan, et \'evitant $n-1$ hyperplans en position g\'en\'erale que l'on peut prendre comme les hyperplans de coordonn\'ees interm\'ediaires. Donc $u_2,...,u_n$ sont harmoniques. Enfin on choisit un autre hyperplan en position g\'en\'erale par rapport \`a ces $n-1$ hyperplans, par exemple le dernier hyperplan de coordonn\'ees. Sa fonction de comptage est $N=\oint_r u_{n+1}$. On veut comparer $N$ et $T$. On raisonne par r\'ecurrence. Posons $u_k^*=\text{Max}(u_k,...,u_{n+1})$. Noter que $T=\oint_ru_1^*$ et $N=\oint_ru_{n+1}^*$. L'hypoth\`ese de r\'ecurrence se compose de l'estim\'ee de $T$ en fonction de $u_k^*$ et d'estim\'ees locales dans l'esprit du lemme 2.

\null\noindent
$(H_k)$ \hskip.5cm 1) $T=\oint_ru_k^*+O(r)$.

\hskip.9cm 2) Pour tout $z$, ou bien $u_k^*(z)=u_1^*(z)$, ou bien il existe $i<k$ et un disque $d$ centr\'e en $z$ de rayon $\rho=O(\vert z\vert)$ tel que $2d\subset (u_k^*<u_i)$ et $\oint_{\gamma}(u_i-u_k^*)=O(\rho)$ pour tout arc $\gamma$ de $d$ de rayon $\geq \rho$.

\null
La conclusion voulue est l'estim\'ee 1) de $(H_{n+1})$. $(H_1)$ est tautologique. Montrons que $(H_{k})$ implique $(H_{k+1})$. Commen\c cons par le 2).

\null
Supposons d'abord que $u_{k+1}^*(z)\geq u_k(z)$. Si $u_{k+1}^*(z)< u_1^*(z)$ alors $u_k^*(z)< u_1^*(z)$ et par hypoth\`ese il existe $i<k$ et $d$ centr\'e en $z$ de rayon $\rho$ avec $2d\subset(u_k ^*<u_i)$ et $\oint_{\gamma}(u_i-u_k^*)=O(\rho)$. Par le lemme d'oubli (lemme 4) on a une conclusion du m\^eme type en rempla\c cant $u_k^*$ par $u_{k+1}^*$, ce qui conclut.

\null
Sinon on a $u_{k+1}^*(z)<u_k(z)$. Soit $2d$ le disque maximal centr\'e en $z$ contenu dans $(u_{k+1}^*<u_k)$. Le rayon $\rho$ de $d$ est contr\^ol\'e par $\vert z \vert$ puisque le lieu $(u_{k+1}^*=u_k)$ est non vide (sans quoi $f$ serait d\'eg\'en\'er\'ee). Au bord de $2d$ on a un point $z'$ dans $(u_{k+1}^*=u_k)$.

\null
Supposons que $u_k^*(z')=u_1^*(z')$. On v\'erifie qu'alors $u_{k+1}^*$ est diff\'erentiable en $z'$ et $\vert \nabla u_{k+1}^*-\nabla u_k\vert(z')$ est contr\^ol\'e par la constante de Brody initiale. Autrement dit $(u_1,...,u_k,u_{k+1}^*)$ se comporte comme les coordonn\'ees logarithmiques d'une courbe de Brody de $P^k(\C)$. Par le lemme 2 appliqu\'e \`a $u_k-u_{k+1}^*$ on obtient l'estim\'ee voulue.

\null
Sinon on a $u_k^*(z')<u_1^*(z')$. Par hypoth\`ese il existe $i<k$ et $d'$ centr\'e en $z'$ de rayon $\rho'$ avec $2d'\subset(u_k^*<u_i)$ et $\oint_{\gamma'}(u_i-u_k^*)=O(\rho')$ pour tout arc $\gamma'$ de $d'$ de rayon $\geq \rho'$. Par le lemme d'oubli (lemme 4) on a la m\^eme conclusion en rempla\c cant $u_k^*$ par $u_{k+1}^*$ et avec $\gamma'$ de rayon $\geq \frac {\rho'}2$. A fortiori $\oint_{\gamma'}(u_k-u_{k+1}^*)^+=O(\rho')$. Maintenant ou bien $d'$ ne contient pas $2d$ ($\rho'<4\rho$) et l'estim\'ee pr\'ec\'edente se transf\`ere \`a $d$. En prenant $\gamma'=2d \cap \bo \frac{d'}2$ dans le lemme 3 on obtient $\oint_{\gamma}(u_k-u_{k+1}^*)=O(\rho)$ . Ou bien $d'$ contient $2d$. Le disque $d'$ contient alors un disque $d''$ centr\'e en $z$ de rayon $\rho'' \geq \frac{\rho'}2$. On a directement par restriction $\oint_{\gamma''}(u_i-u_{k+1}^*)=O(\rho'')$ pour tout arc $\gamma''$ de $d''$ de rayon $\geq \rho''$ (et donc $\oint_{\gamma''}(u_k-u_{k+1}^*)^+=O(\rho''))$.

\null
Montrons maintenant le 1). Remarquons que l'on vient de voir l'alternative suivante pour un point $z$ de $\gamma_r$. Ou bien $u_{k+1}^*(z)\geq u_k(z)$, ou bien il existe un disque $d$ centr\'e en $z$ de rayon $\rho$ tel que $\oint_{\gamma_r\cap d}(u_k-u_{k+1}^*)^+=O(\rho)$. Par le lemme de recouvrement (voir plus haut) on obtient que  $\oint_r(u_k-u_{k+1}^*)^+=O(r)$. D'o\`u $\oint_ru_k^*=\oint_ru_{k+1}^* +O(r)$ ce qui conclut.

\null\noindent
\subhead 2. Le th\'eor\`eme 3 \endsubhead

\null

Il repose sur une technique de renormalisation due \`a Eremenko (voir [\E],[\ES]). On va montrer un second th\'eor\`eme non int\'egr\'e sous la forme $$(q-n+1)\, a(r)\leq \sum n_{H_i}(r)+o(r^2).$$ Ici $a(r)=\text{aire}(f(D_r))$ et $n_{H_i}(r)=$card$(f^{-1}(H_i)\cap D_r)$ (avec multiplicit\'es). Cela entra\^ \i ne directement le th\'eor\`eme 3 par int\'egration.

\null
Soit $(l_i=0)$ l'\'equation de l'hyperplan $H_i$. On pose $u_i= \log \vert l_i \circ f\vert$ et $u=\frac12\log \sum \vert f_i\vert^2$ (o\`u $f=(f_1:...:f_{n+1})$ comme plus haut). Noter que $\Delta u$ est born\'e puisque $f$ est Brody. On veut voir que $$\limsup \frac1{r^2}[(q-n+1)\, \Delta u(D_r)-\sum \Delta u_i(D_r)]\leq 0.$$ 

Autrement dit soit $(r_k)$ une suite de rayons telle que $\frac1{r_k^2}[(q-n+1)\, \Delta u(D_{r_k})-\sum \Delta u_i(D_{r_k})]$ converge. Il faut voir que sa limite est n\'egative ou nulle. Pour cela on note $\lambda_k$ l'homoth\'etie du plan de rapport $r_k$ et on regarde la suite de fonctions $\frac1{r_k^2}(u-h_k)\circ\lambda_k$ o\`u $h_k$ est la plus petite majorante harmonique de $u$ sur $D_{2r_k}$. On obtient des fonctions sousharmoniques $\leq 0$ sur $2D$. Par compacit\'e (voir [\H]), ou bien cette suite tend localement uniform\'ement vers $-\infty$, ou bien quitte \`a extraire elle converge dans $L^1_{\text{loc}}$ vers une fonction sousharmonique $v$. Le premier cas est exclu ici. En effet $h_k(0)=\oint_{2r_k}u=O(r_k^2)$ puisque $f$ est Brody, donc $\frac1{r_k^2}(u(0)-h_k(0))$ reste born\'e. On proc\`ede de m\^eme pour les fonctions $\frac1{r_k^2}(u_i-h_k)\circ\lambda_k$ et on note $v_i$ la limite correspondante. 

\null
Pour $I\subset \{1,...,q\}$ on note $u_I=$ Max$_{\,I}\, u_i$ et
$v_I=$ Max$_{\,I}\, v_i$. Si card$(I)\geq n+1$ on a $u=u_I +O(1)$ par position g\'en\'erale, donc $v=v_I$. En fait $v=v_I$ d\`es que card$(I)=n$.
En effet soit $\gamma$ un cercle dans $2D$. Par le m\^eme argument que pour le th\'eor\`eme 1 (voir aussi l'hypoth\`ese de r\'ecurrence $(H_2)$) on a
$\oint_{\lambda_k(\gamma)}u_{I'}=\oint_{\lambda_k(\gamma)}u_I + O(r_k)$ quand $I\subset I'\subset \{1,...,q\}$ avec card$(I)=n$ et card$(I')=n+1$. \`A la limite $\oint_{\gamma}v_{I'}=\oint_{\gamma}v_I$. Ceci \'etant valable pour tout cercle, en faisant tendre leur rayon vers z\'ero on obtient bien $v_{I'}=v_I$.
Autrement dit en tout point de $2D$ au moins $q-n+1$ fonctions $v_i$ co\" \i ncident avec $v$. 

Pour conclure on va appliquer le lemme de potentiel suivant.
\th{Lemme 5} Soit $w$ une fonction $\delta$-sousharmonique sci $\geq0$ sur $2D$ telle que $\Delta w$ soit $L^{\infty}$ sur $(w=0)$. 
Alors $\Delta w=0$ sur $(w=0)$.\fth

Une fonction $\delta$-sousharmonique est la diff\'erence de deux fonctions sousharmoniques et {\it sci} signifie semicontinue inf\'erieurement. 

\null
Admettons ce r\'esultat pour l'instant. Soit $w=v-v_i$. Elle est bien $\delta$-sousharmonique $\geq 0$ sur $2D$ et sci. En effet comme $\Delta u$ est born\'ee on a la m\^eme borne uniforme pour le laplacien de $\frac1{r_k^2}(u-h_k)\circ\lambda_k$ donc $\Delta v$ est $L^{\infty}$. Par r\'egularit\'e elliptique du laplacien, $v$ est un peu moins que $C^2$ donc continue et $w$ est sci. De plus $\Delta w$ est $L^{\infty}$ sur $(w=0)$. En effet par le lemme de Grishin ([\G],[\ES]) on a $\Delta w\geq 0$ sur $(w=0)$. Donc, toujours sur $(w=0)$, $\Delta v\geq  \Delta v_i$ et $\Delta v_i$ est $L^{\infty}$. On peut donc appliquer le lemme et on obtient que $\Delta v = \Delta v_i$ sur $(v=v_i)$. Maintenant, puisqu'au moins $q-n+1$ fonctions $v_i$ co\" \i ncident avec $v$ en chaque point, on a $(q-n+1)\, \Delta v \leq \sum \Delta v_i$ partout. 

\null
Traduisons ceci avant la limite. Donnons-nous $\epsilon >0$ et $\chi\geq 0$ une fonction plateau sur $(1-\epsilon)D$ \`a support dans $D$. Pour $k$ assez grand on a $$(q-n+1)\int  \chi \, \Delta (u\circ \lambda_k)\leq \sum  \int   \chi \, \Delta (u_i\circ \lambda_k)  +\epsilon\,  r_k^2.$$

Donc $(q-n+1)\, a(r_k(1-\epsilon))\leq \sum n_{H_i}(r_k)+\epsilon\,  r_k^2$. Comme $f$ est Brody il s'ensuit que
$\frac1{r_k^2}[(q-n+1)\, a(r_k)-\sum n_{H_i}(r_k)] \leq C\, \epsilon$ o\`u $C$ est une constante ne d\'ependant que de $q,n$ et la constante de Brody de $f$. On conclut en faisant tendre $\epsilon$ vers 0.

\null
Terminons par la d\'emonstration du lemme 5. Elle suit celle du lemme 2 dans [\ES]. Par une r\'eduction analogue on se ram\`ene au cas o\`u $(w=0)=\bo (2D)\cup K$ avec $K$ compact dans $\overset \circ \to {2D}$. Dans une composante connexe $U$ du compl\'ementaire de $K$ on a donc $w=\int_U g_t\, \Delta w(t)$. Ici $g_t$ est la fonction de Green sousharmonique sur $U$ avec p\^ole en $t$. On la prolonge par $0$ en une fonction $\delta$-sousharmonique du plan. On a $\Delta \, g_t =\delta_t-\omega_t$ o\`u $\omega_t$ est la mesure harmonique sur $\bo U$ par rapport \`a $t$. En sommant sur les composantes $w=\sum \, \int_U g_t\, \Delta w(t)$ et $\Delta w =\sum \, [{\bold 1}_U\Delta w -\, \int_U\omega_t\, \Delta w(t)]$. Donc ${\bold 1}_K\, \Delta w=-{\bold 1}_K \sum \, \int_U \omega_t\, \Delta w(t)$. Or dans chaque $U$ les mesures $\omega_t$ sont absolument continues entre elles et port\'ees
par un bor\'elien de mesure de Lebesgue nulle [\OK]. C'est donc aussi le cas pour ${\bold 1}_K\, \Delta w$. Mais par hypoth\`ese ${\bold 1}_K\, \Delta w$
est absolument continue par rapport \`a la mesure de Lebesgue. Elle est donc nulle.

\Refs
\widestnumber\no{99}
\refno=0
\bref \by M. Barrett et A. Eremenko\paper A generalization of a theorem of Clunie and Hayman\jour Proc. Amer. Math. Soc.\vol140\yr2012\pages1397--1402
\endref
\bref \by J. Clunie et W. Hayman\paper The spherical derivative of integral and meromorphic functions\jour Comm. Math. Helv.\vol40\yr1966\pages117--148
\endref
\bref \by A. Eremenko\paper Brody curves omitting hyperplanes\jour Ann. Acad. Sci. Fenn. Math. \vol35\yr2010\pages565--570
\endref
\bref \by A. Eremenko et M. Sodin\paper The value distribution of meromorphic functions and meromorphic curves from the point of view of potential theory\jour St. Petersburg Math. J.\vol13\yr1992\pages109--136
\endref
\bref \by A. Grishin\paper Sets of regular growth of entire functions. I (Russian)\jour Teor. Funktsiĭ Funktsional. Anal. i Prilozhen.\vol40\yr1983\pages36--47
\endref
\bref \by L. H\"ormander \book The analysis of linear partial differential operators I\publ Springer \yr 1990 \publaddr Berlin
\endref
\bref \by S. Kobayashi \book Hyperbolic complex spaces \bookinfo 
 Grund. Math. Wiss.\vol318 \publ 
 Springer \yr 1998 \publaddr Berlin 
 \endref
\bref \by P. Mattila \book Geometry of sets and measures in euclidean spaces\publ Cambridge University Press \yr 1995 \publaddr Cambridge
\endref
\bref \by R. Nevanlinna \book Analytic functions \bookinfo Grund. Math. Wiss. \vol162 \publ Springer \yr 1970 \publaddr Berlin
\endref
\bref \by B. \O ksendal\paper Null sets for measures orthogonal to $R(X)$\jour Amer. J. Math.\vol94\yr1972\pages331--342
\endref
\bref \by T. Ransford \book Potential theory in the complex plane \publ 
 Cambridge University Press \yr 1995 \publaddr Cambridge 
 \endref

\endRefs

\enddocument